# ON THE SILHOUETTE OF BINARY SEARCH TREES

By Rudolf Grübel

*Leibniz Universität Hannover*

A zero-one sequence describes a path through a rooted directed binary tree $T$; it also encodes a real number in $[0, 1]$. We regard the level of the external node of $T$ along the path as a function on the unit interval, the *silhouette* of $T$. We investigate the asymptotic behavior of the resulting stochastic processes for sequences of trees that are generated by the binary search tree algorithm.

**1. Introduction.** Let $(\xi_n)_{n \in \mathbb{N}}$ be a sequence of independent random variables, where each $\xi_n$ is uniformly distributed on the unit interval. The binary search tree (BST) algorithm sequentially stores these variables in a sequence $(T_n)_{n \in \mathbb{N}}$ of rooted, directed, labeled binary trees. $T_1$ consists of the root node only, with label $\xi_1$. In order to obtain $T_{n+1}$ from $T_n$, we compare $\xi_{n+1}$ with the labels of the nodes along a path through $T_n$, beginning at the root and moving to the left if $\xi_{n+1}$ is smaller, to the right if it is greater than the label associated with the respective node. Once an empty node has been found, we attach it to the tree and $\xi_{n+1}$ is the label of the new node (formal definitions will be given below). The BST algorithm is one of the basic and classical search procedures and is discussed in the standard texts in this area; see, for example, Knuth (1973), Mahmoud (1992) and Sedgewick and Flajolet (1996).

Let $\mathcal{T}_n$ be the set of rooted directed binary trees with $n$ nodes. Then, $T_n$ is a random variable with values in $\mathcal{T}_n$, but the distribution of $T_n$ is *not* the uniform distribution on $\mathcal{T}_n$. For uniformly distributed plane trees or, more generally, simply generated trees, there are various codings, for example, by depth-first search, that relate the trees to random walks (Harris correspondence). These codings provide the basis for an in-depth study of simply generated trees, leading to limit results that involve Brownian excursions, and that, in turn, have applications to certain nonlinear partial differential equations; see Aldous (1991), Chapter 6 in Pitman (2006) and Le Gall (1999).

For the present BST case, we investigate an encoding by a function that we call the *silhouette* $X(T) = (X_s(T))_{0 \leq s \leq 1}$ of the tree $T$. Any $s \in [0, 1]$ defines a path through $T$ by its binary expansion and $X_s(T)$ is simply the depth (or level) of the external node of $T$ along this path. In contrast to other notions such as the profile of the tree [see Chauvin, Drmota and Jabbour-Hattab (2001)], this is a coding in the sense that $T$ can be reconstructed from $X$. The notion and use of paths through the tree is of course not new and appears in, for example, Pittel (1985, 1986); the label "silhouette" was coined in Grübel (2005). Applied to the output sequence $T_n$, $n \in \mathbb{N}$, of the BST algorithm, the silhouette yields a sequence $(X_s(T_n))_{0 \leq s \leq 1}$, $n \in \mathbb{N}$, of stochastic processes. Our main result shows that these processes converge in a weak sense to a nondegenerate limit process as $n \to \infty$.

Section 2 contains some formal definitions related to trees. In Section 3, we define the silhouette and show that distributional convergence to a nondegenerate limit process does not hold with respect to pointwise convergence on the underlying function space. The weak convergence essentially refers to the integrated silhouette, and a key role for the analysis of the latter is played by the discounted external path length, which we discuss in Section 4. In Section 5, we consider the finite-dimensional distributions of the integrated silhouette, prove the convergence to a limit process, characterize the limit distribution as the unique solution to a fixed point equation on a suitable space of measures and study the paths of the limit process. In the final section, we collect some comments on related questions and on possible variations of our findings.

Throughout, we write $\#A$ for the number of elements of the set $A$ and $\mathcal{L}(X)$ for the distribution of the random variable $X$. Sometimes, we write $X \sim \mu$ instead of $\mathcal{L}(X) = \mu$. A random variable $X$ is stochastically smaller than or equal to another random variable $Y$ (both real-valued), written $X \leq_{\mathcal{D}} Y$, if $P(X \geq x) \leq P(Y \geq x)$ for all $x \in \mathbb{R}$. Finally, "$=_{\mathcal{D}}$" denotes equality in distribution and "$\to_{\mathcal{D}}$" denotes convergence in distribution.

**2. Some notation for trees.** A tree is a graph and thus consists of vertices (or nodes) and edges. In the context of binary trees, it is convenient to represent (or define) nodes as elements of $\mathcal{N}$,

$$\mathcal{N} := \bigcup_{k=0}^{\infty} \{0, 1\}^k,$$

where $\{0, 1\}^0 := \{\varnothing\}$. Stated in different terminology, $\mathcal{N}$ is the set of finite words over the alphabet that consists of the two letters 0 and 1. By a rooted, directed binary tree, we mean a finite set $T$ of (internal) nodes with the following property:

$$(1) \qquad u = (u_1, \ldots, u_k) \in T, \qquad k \geq 1 \quad \Longrightarrow \quad \tilde{u} := (u_1, \ldots, u_{k-1}) \in T.$$

A binary tree can therefore be regarded as a finite and prefix-stable set of finite words with letters 0 and 1. Informally, $u_i = 0$ means that we move to the left (to the right, if $u_i = 1$) on the path from the root node to the node in question.

We may interpret $\tilde{u}$ in (1) as the direct ancestor or predecessor of $u$. The root node is represented by the empty string and has no predecessor. The edges of the tree are the pairs $(\tilde{u}, u)$ with $u \neq \varnothing$. The size of a tree is simply the number of its nodes. $\mathcal{T}$ denotes the set of all binary trees. By a labeled tree, we mean a pair $(T, \phi)$, with $T \in \mathcal{T}$ and a function $\phi \colon T \to \mathbb{R}$; the value $\phi(u)$ is the label associated with the node $u \in T$.

Given a tree $T \in \mathcal{T}$, we may now formally define two associated trees, $L(T)$ and $R(T)$, the left and right subtree of $T$, by

$$(2) \qquad \begin{aligned} L(T) &:= \{u = (u_1, \ldots, u_k) \in \mathcal{N} : (0, u_1, \ldots, u_k) \in T\}, \\ R(T) &:= \{u = (u_1, \ldots, u_k) \in \mathcal{N} : (1, u_1, \ldots, u_k) \in T\}. \end{aligned}$$

Obviously, any nonempty $T \in \mathcal{T}$ is uniquely determined by the corresponding subtrees $L(T)$ and $R(T)$ (which may, of course, be empty). For $u \in \mathcal{N}$ and $T \in \mathcal{T}$, let $T(u)$ be the subtree of $T$ that consists of $u$, now regarded as the root node, and all descendants of $u$ in $T$. A formal definition, as in (2), is straightforward. Indeed, we may consider $L(T)$ and $R(T)$ as the subtrees $T(u)$ associated with the nodes $u = (0)$ and $u = (1)$, respectively.

A node $u = (u_1, \ldots, u_k)$ has depth $|u| = k$. It is an external node of $T$ if $u$ itself is not an element of $T$, but its predecessor is. The formal definition of the set $\partial T$ of external nodes of the tree $T \in \mathcal{T}$ is

$$\partial T := \{u = (u_1, \ldots, u_k) \in \mathcal{N} : k \geq 1, u \notin T, (u_1, \ldots, u_{k-1}) \in T\}.$$

We augment this with the convention that $\partial T_0 = \{\varnothing\}$ for the empty tree $T_0$.

**3. The silhouette.** A sequence $u = (u_k)_{k \in \mathbb{N}}$ in $\{0, 1\}^{\mathbb{N}}$ can be regarded as the binary expansion $s = \sum_{k=1}^{\infty} u_k 2^{-k}$ of some number $s$ in the unit interval $[0, 1]$. Conversely, for any $s \in [0, 1]$, we have a unique such binary expansion if we require that binary rationals $s < 1$ have $u_k = 0$ for all $k \geq k_0$, for some $k_0 \in \mathbb{N}$. For any $T \in \mathcal{T}$, we now introduce its *silhouette* as the function $s \mapsto X_s(T)$ on the unit interval defined by

$$X_s(T) := \min\{k \in \mathbb{N}_0 : (u_1, \ldots, u_k) \notin T\};$$

an informal description of the silhouette was given in Section 1. These functions are piecewise constant on intervals with binary rational endpoints if the length of these intervals is chosen small enough, and they are continuous from the right and have left-hand limits; at $s = 1$ they are left-continuous. The left part of Figure 1 shows an example for a tree of size 500.

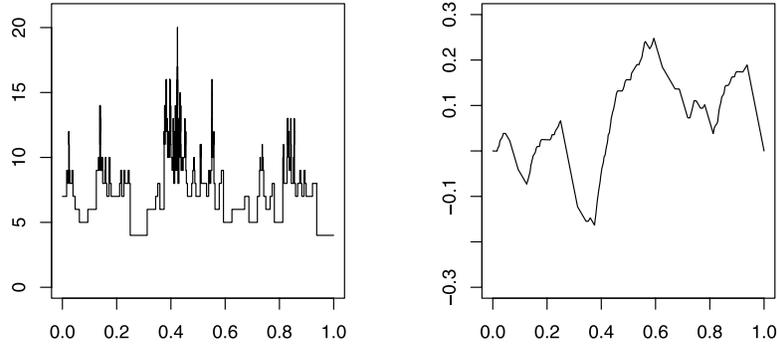

Fig. 1.   *An example of a silhouette (left) and the corresponding normalized integrated silhouette (right).*

We will occasionally write $X(T, s)$ instead of $X_s(T)$. A basic fact is the following recursion, which, because of the preparations undertaken in the previous section, can now be expressed concisely as

$$(3) \qquad X_s(T) = 1 + X_{2s}(L(T))1_{[0,1/2)}(s) + X_{2s-1}(R(T))1_{[1/2,1]}(s)$$

for $0 \leq s \leq 1$, provided that $T \neq \varnothing$. Obviously, $X(\varnothing) \equiv 0$.

Now, let $T_n$ be the random tree generated by the BST algorithm from the first $n$ variables in a sequence of independent, uniformly distributed random variables, as explained in the Introduction [hereafter, we will simply refer to $(T_n)_{n\in\mathbb{N}}$ as a *BST sequence*]. Then, for each $n \in \mathbb{N}$, the silhouette $X(T_n) = (X_s(T_n))_{0 \leq s \leq 1}$ can be regarded as a stochastic process with time parameter ranging over the unit interval.

We first consider the finite-dimensional distributions of the silhouette processes. It is well known [see, e.g., Régnier (1989)] that the present combination of input and algorithm leads to the following stochastic dynamics of the tree sequence: $T_{n+1}$ is obtained from $T_n$ by picking one of the $n + 1$ external nodes of $T_n$ uniformly at random and then adding this node to $T_n$. As a consequence, the sequence $(X_s(T_n))_{n\in\mathbb{N}}$ of depths of the external nodes along the path $s$ can be represented as the sequence of partial sums of an independent sequence $(I_n)_{n\in\mathbb{N}}$ of indicator variables, with $P(I_n = 1) = 1/n$ for all $n \in \mathbb{N}$. Extending this observation to more than one path provides the basis for our first result.

THEOREM 1.   *Let $s_1, \ldots, s_d$, with $0 \leq s_1 < \cdots < s_d \leq 1$, be given. The $d$-dimensional random vectors $Y_n = (Y_{n,1}, \ldots, Y_{n,d})$ with*

$$Y_{n,l} := \frac{X(T_n, s_l) - \log n}{\sqrt{\log n}} \qquad for \; l = 1, \ldots, d$$

*then converge in distribution as $n \to \infty$ to a $d$-dimensional standard normal limit.*

PROOF. Given $s_1, \ldots, s_d$, we define a sequence $(Z_n)_{n \in \mathbb{N}}$ of $d$-dimensional random vectors $Z_n = (Z_{n,1}, \ldots, Z_{n,d})$ by $Z_{n,l} = X(T_n, s_l) - X(T_{n-1}, s_l)$. As we add one node at a time, these all take values in the set consisting of the zero vector and the standard basis vectors $e_j$, $j = 1, \ldots, d$, where the $k$th component of $e_j$ is 1 if $k = j$ and 0 otherwise. Let $(\tilde{Z}_n)_{n \in \mathbb{N}_0}$ be another such sequence where, now, $\tilde{Z}_n = (\tilde{Z}_{n,1}, \ldots, \tilde{Z}_{n,d})$, $n \in \mathbb{N}$, are independent, with $P(\tilde{Z}_n = e_l) = 1/n$ for $l = 1, \ldots, d$ and $P(\tilde{Z}_n = (0, 0, \ldots, 0)) = 1 - l/n$ if $n \geq d$. For the "tilded" variables, we have asymptotic normality by the multivariate version of the Lindeberg–Feller central limit theorem (or the usual one-dimensional version, together with the Cramér–Wold device). The difference between the standardized partial sum processes associated with the $Z$- and the $\tilde{Z}$-variables is asymptotically negligible as the time that the last of the pairwise last common ancestors is reached is finite with probability 1. □

Hence, if we consider the silhouette itself, the appropriate scaling would lead to independent components in the limit. This shows that in order to obtain a nondegenerate limit process for the silhouette sequence $(X(T_n))_{n \in \mathbb{N}}$ associated with a BST sequence $(T_n)_{n \in \mathbb{N}}$, we need to weaken the notion of convergence. A standard strategy is to regard the paths of the process $X(T_n)$ as "weak" functions in the sense of linear forms on some function space $\mathcal{F}$, that is, to investigate the random linear functionals $f \mapsto \int_0^1 X_t(T_n) f(t) \, dt$, $f \in \mathcal{F}$. A key role is played by $f \equiv 1$, a case that we study in the next section.

**4. The discounted external path length.** For a binary tree $T$, let

$$U_k(T) := \#\{u \in \partial T : |u| = k\}$$

be the number of external nodes of $T$ with depth $k$, $k \in \mathbb{N}_0$. We then have

$$\eta(T) := \int_0^1 X_s(T) \, ds = \sum_{k=1}^{\infty} 2^{-k} k U_k(T),$$

which means that we can regard the integral of the silhouette over the whole unit interval as a discounted external path length of the tree.

Now, let $(T_n)_{n \in \mathbb{N}}$ be a BST sequence, as explained in the previous section; we abbreviate $\eta(T_n)$ to $\eta_n$. For the proof of Theorem 1, we have used the dynamical view of this sequence, obtaining $T_{n+1}$ from $T_n$ by inclusion of a randomly chosen element of $\partial T_n$. For the analysis of the BST sequence, the recursive view is equally important: for $n \geq 1$, the subtrees $L(T_n)$ and $R(T_n)$ are conditionally independent given $I_n := \#L(T_n)$, $I_n$ is uniformly distributed on $\{0, \ldots, n-1\}$ and, on $I_n = k$, $L(T_n) =_{\mathcal{D}} T_k$, $R(T_n) =_{\mathcal{D}} T_{n-1-k}$. This is a consequence of the fact that the subsequences of the input sequence $(\xi_n)_{n \in \mathbb{N}}$ to the BST algorithm that consist of the values smaller than and greater than $\xi_1$ are independent conditionally on $\xi_1$ and uniformly distributed on the intervals $(0, \xi_1)$ and $(\xi_1, 1)$, respectively. Hence, we obtain from (3) (or directly) that, for $n \geq 1$,

$$(4) \qquad \eta_n =_{\mathcal{D}} 1 + \tfrac{1}{2}(\eta_{I_n} + \eta'_{n-1-I_n}),$$

with $(\eta_m)_{m \in \mathbb{N}_0}$, $(\eta'_m)_{m \in \mathbb{N}_0}$ and $I_n$ independent, $\eta'_m =_{\mathcal{D}} \eta_m$ for all $m \in \mathbb{N}_0$ and $I_n \sim \mathrm{unif}\{0, 1, \ldots, n-1\}$. Clearly, $\eta_0 \equiv 0$. Let $H(n) = H_n = \sum_{k=1}^{n} 1/k$ be the $n$th harmonic number. For $a_n := E\eta_n$, we have $a_0 = 0$ and (4) implies that

$$a_n = 1 + \frac{1}{n} \sum_{m=0}^{n-1} a_m \qquad \text{for all } n \in \mathbb{N},$$

which easily leads to

$$(5) \qquad\qquad\qquad\qquad E\eta_n = H_n.$$

The undiscounted path length $\sum_{k=1}^{\infty} k U_k(T_n)$ plays a key role in the analysis of Quicksort and we now use the techniques that were successful in that situation: martingales [see Régnier (1989)] and the contraction method [see Rösler (1991)]. The filtration $(\mathcal{F}_n)_{n \in \mathbb{N}}$ of interest in the former context will be the one generated by the sequence $(T_n)_{n \in \mathbb{N}}$. We write $L^2$ for the set of square integrable random variables.

THEOREM 2. *As $n \to \infty$, $\eta_n - H_n$ converges almost surely and in quadratic mean to a random variable $\eta_\infty$. Within the set of distributions with finite second moment and zero mean, the distribution of $\eta_\infty$ is characterized by the fixed point equation*

$$(6) \qquad \eta_\infty =_{\mathcal{D}} \tfrac{1}{2}(\eta_\infty + \eta'_\infty) + \zeta_\infty,$$

*where $\eta_\infty$, $\eta'_\infty$ and $\zeta_\infty$ are independent, $\eta_\infty =_{\mathcal{D}} \eta'_\infty$ and*

$$(7) \qquad \zeta_\infty := 1 + \tfrac{1}{2}(\log(\xi) + \log(1 - \xi)),$$

*with $\xi$ uniformly distributed on the unit interval.*

PROOF. The transition from $T_n$ to $T_{n+1}$ means that an external node of (random) level $K$ becomes an internal node. This entails a loss of $K 2^{-K}$, but, as the new internal node spawns two external nodes at level $K+1$, there is also a gain of $2(K+1)2^{-K-1}$ for the discounted external path length, hence

$$\eta_{n+1} - \eta_n = 2(K+1)2^{-K-1} - K 2^{-K} = 2^{-K}.$$

By the stochastic dynamics of the tree sequence described in Section 3, we have that, given $T_n$ with the associated values $U_k(T_n)$ for the number of external nodes at level $k$,

$$P[K = k|T_n] = \frac{1}{n+1} U_k(T_n) \qquad \text{for all } k \in \mathbb{N}_0.$$

Hence, with $\mathcal{F}_n$ as defined above,

$$(8) \qquad E[\eta_{n+1} - \eta_n | \mathcal{F}_n] = \frac{1}{n+1} \sum_{k=0}^{\infty} U_k(T_n) 2^{-k} = \frac{1}{n+1}.$$

[The last equality uses the well-known fact that $\sum_{k=0}^{\infty} 2^{-k} U_k(T) = 1$ for all binary trees $T$. Note that $2^{-k} U_k(T_n)$ is the Lebesgue measure of $X(T_n)^{-1}(\{k\})$, so this fact has a simple interpretation in terms of the silhouette.] From (8), we immediately obtain that $(\eta_n - H_n, \mathcal{F}_n)_{n \in \mathbb{N}_0}$ is a zero mean martingale; (8) also provides an alternative proof for (5). The individual random variables are all bounded and hence elements of $L^2$.

We next show that the martingale is bounded in $L^2$. Let $\sigma_n^2 := E(\eta_n - H_n)^2 = \text{var}(\eta_n)$. From (4), we obtain

$$\sigma_n^2 = \tfrac{1}{4} \text{var}(\eta_{I_n} + \eta'_{n-1-I_n}).$$

Because of

$$E(\text{var}[\eta_{I_n} + \eta'_{n-1-I_n} | I_n]) = E(\sigma_{I_n}^2 + \sigma_{n-1-I_n}^2) = \frac{2}{n} \sum_{m=0}^{n-1} \sigma_m^2$$

and

$$\text{var}(E[\eta_{I_n} + \eta'_{n-1-I_n} | I_n]) = \text{var}(H_{I_n} + H_{n-1-I_n}),$$

the conditional variance formula leads to

$$\sigma_n^2 = \frac{1}{2n} \sum_{m=0}^{n-1} \sigma_m^2 + \tfrac{1}{4} b_n, \qquad \text{with } b_n := \text{var}(H_{I_n} + H_{n-1-I_n}).$$

As $((\eta_n - H_n)^2, \mathcal{F}_n)_{n \in \mathbb{N}_0}$ is a submartingale, we have that $m \mapsto \sigma_m^2$ is non-decreasing, so the sum may be bounded from above by $\sigma_{n-1}^2/2$. To obtain boundedness of the sequence $(\sigma_n^2)_{n \in \mathbb{N}_0}$, it is therefore enough to show that the sequence $(b_n)_{n \in \mathbb{N}}$ is bounded. This, in turn, will follow from the boundedness of $(E(H_{I_n} - H_n)^2)_{n \in \mathbb{N}}$ if we use Minkowski's inequality and the fact that $I_n$ and $n - 1 - I_n$ have the same distribution. From the elementary inequalities

$$\log m \leq H_m \leq \log m + 1 \qquad \text{for all } m \in \mathbb{N},$$

we obtain

$$\log \frac{I_n}{n} - 1 \leq H_{I_n} - H_n \leq \log \frac{I_n}{n} + 1$$

on $\{I_n > 0\}$. Hence, the required boundedness will be implied by

$$\sup_{n \in \mathbb{N}} E\left(1_{\{I_n > 0\}} \log \frac{I_n}{n}\right)^2 < \infty,$$

which finally follows from

$$E\left(1_{\{I_n > 0\}} \log \frac{I_n}{n}\right)^2 = \frac{1}{n}\sum_{m=1}^{n-1}\left(\log \frac{m}{n}\right)^2 \leq \int_0^1 (\log x)^2 \, dx = 2.$$

Because of its boundedness in $L^2$, the martingale $(\eta_n - H_n)_{n \in \mathbb{N}}$ converges to some limit $\eta_\infty$ almost surely and in $L^2$ as $n \to \infty$.

Next, we derive the fixed point equation for the distribution of $\eta_\infty$. We first note that (4) implies that for $n \geq 1$,

(9)    $$\eta_n - H_n =_{\mathcal{D}} \tfrac{1}{2}(\eta_{I_n} - H_{I_n}) + \tfrac{1}{2}(\eta'_{n-1-I_n} - H_{n-1-I_n}) + \zeta_n$$

with

$$\zeta_n := 1 + \tfrac{1}{2}(H_{I_n} + H_{n-1-I_n}) - H_n.$$

We may assume that $I_n = \lfloor n\xi \rfloor$ with $\xi \sim \mathrm{unif}(0,1)$. With the standard asymptotic result for harmonic numbers, $H_n = \log(n) + \gamma + o(1)$, we then obtain

$$\zeta_n = 1 + \frac{1}{2}(\log(\lfloor n\xi \rfloor) - \log(n)) + \frac{1}{2}(\log(n-1-\lfloor n\xi \rfloor) - \log(n)) + o(1)$$

$$= 1 + \frac{1}{2}\log\left(\frac{\lfloor n\xi \rfloor}{n}\right) + \frac{1}{2}\log\left(1 - \frac{\lfloor n\xi \rfloor + 1}{n}\right) + o(1)$$

$$\to 1 + \frac{1}{2}\log(\xi) + \frac{1}{2}\log(1-\xi).$$

Using $\eta_n - H_n \to_{\mathcal{D}} \eta_\infty$ and the distributional assumptions on $(\eta_n)_{n \in \mathbb{N}}$, $(\eta'_n)_{n \in \mathbb{N}}$ and $I_n$, we now obtain (6) by letting $n \to \infty$ in (9).

Finally, the right-hand side of (6) defines an operator $\Psi$ on the space $\mathcal{M}_{2,0}$ of distributions $\mu$ with mean zero and finite second moment. A straightforward calculation shows that

$$d_2^2(\Psi(\mu_1), \Psi(\mu_2)) \leq \tfrac{1}{2}d_2^2(\mu_1, \mu_2) \qquad \text{for all } \mu_1, \mu_2 \in \mathcal{M}_{2,0},$$

where $d_2$ denotes the usual Wasserstein 2-distance. Hence, $\Psi$ is a contraction and the distribution of $\eta_\infty$ is characterized in $\mathcal{M}_{2,0}$ by (6). $\square$

Note that with $\xi$ as in the theorem, $-\log(\xi)$ and $-\log(1-\xi)$ both have an exponential distribution with mean 1 (of course, they are not independent); in particular, $\zeta_\infty$ has finite moments of all orders, $E\zeta_\infty = 0$ and $\mathrm{var}(\zeta_\infty) = 1 - \pi^2/12 \approx 0.1775$. The fact that the distribution of $\zeta_\infty$ is nondegenerate implies the same for the distribution of $\eta_\infty$.

It is known that in order to obtain a nondegenerate distributional limit for the undiscounted external path length $\tilde{\eta}_n := \sum_{k=1}^{\infty} k U_{n,k}$, we need to shift and rescale: $(\tilde{\eta}_n - 2n \log n)/n \to_{\mathcal{D}} \tilde{\eta}_\infty$, where $\tilde{\eta}_\infty$ is a real, nonconstant random variable [see Régnier ([1989](#)) and Rösler ([1991](#))]. It may seem surprising that for the discounted version, it is enough to shift. Roughly, because of the asymptotic independence of the marginal distributions of the silhouette, the "thin spikes" and the "broad valleys" cancel out to a certain extent; see Figure [1](#). Note that $\tilde{\eta}_n/(n+1)$ and $\eta_n$ can be regarded as the mean associated with the random distribution $\tilde{\mu}_n$ and $\mu_n$, respectively, where $\tilde{\mu}_n$ has density $k \mapsto 2^k/(n+1)$ with respect to $\mu_n$. In fact, if we regard $\tilde{\eta}_n/(n+1)$ as the basic variables, then, again, it is enough to shift to obtain a nondegenerate limit distribution.

Despite the close connection between $\tilde{\mu}_n$ and $\mu_n$, the fixed point equations for the respective distributional limits are quite different in the two cases. Indeed, in the undiscounted case, Rösler ([1991](#)) obtained the equation

$$(10) \qquad \tilde{\eta}_\infty =_{\mathcal{D}} \xi \tilde{\eta}_\infty + (1-\xi) \tilde{\eta}'_\infty + C(\xi),$$

where $\xi$ is uniformly distributed on the unit interval, $\tilde{\eta}_\infty$, $\tilde{\eta}'_\infty$ and $\xi$ are independent, $\tilde{\eta}_\infty =_{\mathcal{D}} \tilde{\eta}'_\infty$ and

$$(11) \qquad C(x) = 1 + 2(x \log(x) + (1-x) \log(1-x)).$$

A first major difference between ([10](#)) and ([6](#)) is the fact that in the latter case, the linear combination on the right-hand side of the independent copies of the left-hand side has the deterministic coefficients $1/2$ and $1/2$ instead of $\xi$ and $1-\xi$. As we will see in the next result, this makes it possible to obtain a simple and explicit representation for $\eta_\infty$, which is lacking in the undiscounted case. We mention, in passing, that the distribution of $\tilde{\eta}_\infty$ has been the subject of considerable attention; see, for example, Cramer ([1996](#)) Devroye, Fill and Neininger ([2000](#)) and Fill and Janson ([2000](#)). A second important difference is the fact that the function $C$ in ([11](#)) is bounded, whereas the corresponding function of $\xi$ in ([7](#)) is only bounded from above. This has an important consequence for the tail behavior (finiteness domain of the moment generating function) of the respective solutions.

THEOREM 3. *Let $\{\zeta_{n,k} : n \in \mathbb{N}_0, k \in \{1, \ldots, 2^n\}\}$ be a family of independent random variables, all with the same distribution as $\zeta_\infty$, where $\zeta_\infty$ is given in ([7](#)). Then, with $\eta_\infty$ as in Theorem [2](#),*

$$(12) \qquad \eta_\infty =_{\mathcal{D}} \sum_{n=0}^{\infty} 2^{-n} \sum_{k=1}^{2^n} \zeta_{n,k}.$$

*Further, the moment generating function $M(t) = E \exp(t \eta_\infty)$ for $\eta_\infty$ is finite for all $t > -2$. Finally, with $\eta_n$ as in Theorem [2](#), we have that*

$$(13) \qquad M_n(t) := E \exp(t(\eta_n - H_n)) \leq M(t) \qquad \text{for all } t > -2.$$

PROOF. It is easily checked that the sequence $(\eta_{\infty,n})_{n\in\mathbb{N}}$ of partial sums,

$$\eta_{\infty,n} := \sum_{m=0}^{n} 2^{-m} \sum_{k=1}^{2^m} \zeta_{m,k},$$

is a zero mean martingale that is bounded in $L^2$ and hence converges almost surely and in $L^2$ to a random variable $\eta_{\infty,\infty}$. It is equally easy to check that $\eta_{\infty,\infty}$ solves (6), hence (12) follows with the uniqueness statement in Theorem 2.

We know from the proof of Theorem 2 that $(\eta_n - H_n)_{n\in\mathbb{N}}$ is a martingale that converges in $L^2$ to $\eta_\infty$. This convergence implies the representation

(14)           $\eta_n - H_n = E[\eta_\infty \mid \mathcal{F}_n]$      for all $n \in \mathbb{N}$,

where $(\mathcal{F}_n)_{n\in\mathbb{N}}$ denotes the natural filtration associated with $(\eta_n - H_n)_{n\in\mathbb{N}}$. We now know that $\eta_\infty$ and $\eta_{\infty,\infty}$ have the same moment generating function, hence (13) follows from (14) on using Jensen's inequality for conditional expectations.

It remains to prove the finiteness of $M(t)$ for $t > -2$. The moment generating function for $\zeta_\infty$ can be given explicitly as

$$M_0(t) = e^t \int_0^1 x^{t/2}(1-x)^{t/2}\,dx = \frac{\Gamma(1+t/2)^2 e^t}{\Gamma(2+t)},$$

where, of course, $t/2 > -1$ is required for the integral to be finite. From (12), we obtain

$$M(t) = \prod_{n=0}^{\infty} M_0(2^{-n}t)^{2^n}.$$

As $E\zeta_\infty = 0$ and hence $M_0'(0) = 0$, it is straightforward to show that the product converges for $t > -2$.  □

In the undiscounted case, the moment generating function exists on the whole real line. The inequality (13) will be used to obtain a uniform tail bound for the variables $\eta_n - H_n$ that is needed in the next section.

**5. The integrated silhouette.** We now return to the functional point of view explained at the end of Section 3, taking for $\mathcal{F}$ the set of indicator functions $1_{[0,t]}$, $0 \le t \le 1$. This leads us to consider the integrated silhouette $Y(T) = (Y_t(T))_{0\le t\le 1}$ associated with a binary tree $T$, where

$$Y_t(T) := \int_0^t X_s(T)\,ds \qquad \text{for all } t \in [0,1].$$

We define the normalized version $Y^\circ(T) = (Y_t^\circ(T))_{0\le t\le 1}$ by

$$Y_t^\circ(T) := Y_t(T) - tY_1(T) \qquad \text{for all } t \in [0,1].$$

This "ties down" the original process, in the sense that $Y_0^\circ(T) = Y_1^\circ(T) = 0$. Also, $Y_1(T) = \eta(T)$, the discounted external path length discussed in the previous section. We also use $\eta^\circ(T) := \eta(T) - H(\#T)$ to denote the centered external path length.

Suppose, now, that $(T_n)_{n \in \mathbb{N}}$ is a BST sequence. Our aim is a functional limit theorem for the resulting sequence $(Y(T_n))_{n \in \mathbb{N}}$ of stochastic processes. Instead of $Y(T_n)$, we will consider the pairs

$$(15) \qquad Z_n := \begin{pmatrix} Y^\circ(T_n) \\ \eta^\circ(T_n) \end{pmatrix},$$

which we regard as random quantities with values in the linear space

$$(16) \qquad S := C_{00}[0,1] \times \mathbb{R},$$

where $C_{00}[0,1]$ denotes the set of all continuous functions $f : [0,1] \to \mathbb{R}$ that have $f(0) = f(1) = 0$. We can obviously recover $Y(T_n)$ from $Z_n$ (and $H_n$). Together with

$$\|z\| := \|f\|_\infty + |a| \qquad \text{for all } z = (f, a) \in S,$$

the linear space $S$ becomes a separable Banach space. Here, $\|f\|_\infty := \sup_{0 \le t \le 1} |f(t)|$ denotes the supremum norm. We will show that $Z_n$ converges in distribution as $n \to \infty$, where the convergence refers to the topological structure induced by this norm. For this, we follow the classical route, as laid out in Billingsley ([1968](#)), considering first the finite-dimensional distributions and then proving tightness, but the actual details need to be adapted to the present setup.

In connection with the finite-dimensional distributions, instead of considering the $Y^\circ(T_n)$ part of $Z_n$ at arbitrary arguments $t_0, \ldots, t_d \in [0,1]$, we restrict ourselves to complete sets of binary rationals of the same depth, that is, we take $d = 2^k$ and $t_j = j2^{-k}$ for $j = 0, \ldots, 2^k$. A standard argument using the continuity of the paths of $Y^\circ(T_n)$ shows that weak convergence of the resulting random vectors, for all $k \in \mathbb{N}$, is enough to characterize the limit distribution. In fact, in order to simplify the description of the limiting finite-dimensional distributions, we will not consider the values $Y_{t_j}^\circ(T_n)$ themselves, but rather the differences. Again, in the present situation, this suffices as $Y_0^\circ(T_n) \equiv 0$.

Formally, for $k \in \mathbb{N}$, let $\Delta_k$ be the operator that maps a function $f : [0,1] \to \mathbb{R}$ to a vector $\Delta_k f = ((\Delta_k f)_j)_{j=1}^{2^k}$ of dimension $2^k$, with

$$(\Delta_k f)_j := f(j2^{-k}) - f((j-1)2^{-k}) \qquad \text{for } j = 1, \ldots, 2^k.$$

The following theorem shows that these increments converge in distribution and also gives a description of the limits. For the statement of the result, we need some more notation. Let $u(k, j) = (u_1(k, j), \ldots, u_k(k, j)) \in \{0, 1\}^k$, $j =$ $1, \dots, 2^k$, be the nodes of depth $k \in \mathbb{N}$ in the order of their associated binary rationals [the components can be given explicitly as $u_m(k,j) = \lfloor 2^{-k+m}(j-1) \rfloor \pmod 2$] and let

$$u(k,j,l) = (u_1(k,j), \dots, u_l(k,j))$$

be the ancestor of $u(k,j)$ at depth $l$, $l = 0, \dots, k-1$, with the understanding that $u(k,j,0) = \varnothing$. We also write $1_k$ for the $2^k$-dimensional vector that has all entries equal to 1.

THEOREM 4. *For $k \in \mathbb{N}$ fixed and with $n \to \infty$,*

$$(17) \qquad \begin{pmatrix} \Delta_k Y^\circ(T_n) \\ \eta^\circ(T_n) \end{pmatrix} \to_{\mathcal{D}} \begin{pmatrix} 2^{-k}(k1_k + \rho_k + \eta_k - \eta_\infty^\circ) \\ \eta_\infty^\circ \end{pmatrix}.$$

*Here, $\rho_k = (\rho_{k,1}, \dots, \rho_{k,2^k})$ and $\eta_k = (\eta_{k,1}, \dots, \eta_{k,2^k})$ are independent $2^k$-dimensional random vectors, given (distributionally) as follows. The components $\eta_{k,j}$ of $\eta_k$ are independent and each $\eta_{k,j}$ has the same distribution as $\eta_\infty$, defined in (6). Further, the components $\rho_{k,j}$ of $\rho_k$ have the joint distributional representation*

$$(18) \qquad \rho_{k,j} =_{\mathcal{D}} \sum_{\{l\,:\,u_l(k,j)=1\}} \log \xi_{u(k,j,l-1)} + \sum_{\{l\,:\,u_l(k,j)=0\}} \log(1 - \xi_{u(k,j,l-1)}),$$

*where $\xi_u$, $u \in \mathcal{N}$, is a family of independent random variables, all uniformly distributed on the unit interval. Finally,*

$$\eta_\infty^\circ = 2^{-k} \sum_{j=1}^{2^k} (\rho_{k,j} + \eta_{k,j}).$$

PROOF. Let $T^{k,j} := T(u(k,j))$ denote the subtree of $T$ with root $u(k,j)$. We have, provided that the fill level of $T$ is at least $k$ so that none of the subtrees is empty,

$$(19) \qquad \begin{aligned} &Y_{j2^{-k}}^\circ(T) - Y_{(j-1)2^{-k}}^\circ(T) \\ &\quad = \int_{(j-1)2^{-k}}^{j2^{-k}} X_s(T)\,ds - 2^{-k}\eta(T) \\ &\quad = 2^{-k}(k + \eta(T^{k,j}) - \eta(T)) \\ &\quad = 2^{-k}(k + H(\#T^{k,j}) - H(\#T) + \eta^0(T^{k,j}) - \eta^\circ(T)) \end{aligned}$$

for $j = 1, \dots, 2^k$. Further, these differences sum to zero in view of $Y_1^\circ(T) = Y_0^\circ(T)$, hence,

$$(20) \qquad \eta^\circ(T) = 2^{-k} \sum_{j=1}^{2^k} (k + \eta^0(T^{k,j}) + H(\#T^{k,j}) - H(\#T)).$$

The development thus far has been for a fixed tree $T \in \mathcal{T}$. We now substitute the elements $T_n$ of a BST sequence for $T$. Let $N_{n,k} := (N_{n,k,1}, \ldots, N_{n,k,2^k})$ be the random vector that counts the size of the subtrees of $T_n$ at level $k$, that is,

$$N_{n,k,j} = \# T_n^{k,j}, \qquad j = 1, \ldots, 2^k.$$

Because of (19) and (20), in order to obtain the distributional limits in (17), it is enough to work out the asymptotic behavior of the vector

$$(21) \quad (\eta^\circ(T_n^{k,1}), \ldots, \eta^\circ(T_n^{k,2^k}), H(N_{n,k,1}) - H(n), \ldots, H(N_{n,k,2^k}) - H(n))$$

as $n \to \infty$ (the quantities of interest can be written as a fixed linear function of these vectors, so the continuous mapping theorem applies).

Given that a node $u$ at level $l$ has $\# T_n(u) = m$, there are $\lfloor m\xi_u \rfloor$ nodes in the left subtree and $m - 1 - \lfloor m\xi_u \rfloor$ nodes in the right subtree of $T_n(u)$, independent of what happened at levels $0, \ldots, l-1$. Hence,

$$\frac{1}{n} N_{n,k} \to_{\mathcal{D}} V_k = (V_{k,1}, \ldots, V_{k,2^k}),$$

where

$$V_{k,j} =_{\mathcal{D}} \prod_{l=1}^{k} \xi_{u(k,j,l-1)}^{u_l(k,j)} (1 - \xi_{u(k,j,l-1)})^{1-u_l(k,j)},$$

jointly in $j = 1, \ldots, 2^k$. Further, given $N_{n,k,j} = m_j$, the trees $T_n^{k,j}$ are independent with

$$T_n^{k,j} =_{\mathcal{D}} T_{m_j}, \qquad j = 1, \ldots, 2^k,$$

so that the familiar asymptotics of harmonic numbers, together with Theorem 2, imply that the random vector in (21) converges in distribution to

$$(\eta_{k,1}, \ldots, \eta_{k,2^k}, \rho_{k,1}, \ldots, \rho_{k,2^k})$$

as $n \to \infty$, with $\eta_{k,j}, \rho_{k,j}$ as in the statement of the theorem. The theorem now follows by appropriately combining elements. $\quad\square$

Hence, in contrast to the silhouette itself, where the individual random variables are asymptotically independent, we now have limiting finite-dimensional distributions that might be compatible with a limit process that has somewhat regular (e.g., continuous) paths. Below, we will see that the representation of the finite-dimensional distributions given in (18) provides the key to the proof of path properties of the limit process.

To obtain convergence in distribution, we need tightness of the sequence $(Z_n)_{n \in \mathbb{N}}$. For this, we require a technical detail that we state separately as a lemma. We say that a family $(X_i)_{i \in I}$ of nonnegative random variables has *uniformly exponentially decreasing tails* if, for some constants $\kappa > 0$ and $C < \infty$,

$$(22) \qquad P(X_i \geq x) \leq C \exp(-\kappa x) \qquad \text{for all } x \geq 0 \text{ and } i \in I.$$

Some obvious properties of this notion, such as stability with respect to taking sums, will be used below without further comment.

LEMMA 5. *Let $(X_n)_{n \in \mathbb{N}_0}$ and $(A_n)_{n \in \mathbb{N}}$ be sequences of nonnegative random variables with $X_0 \equiv 0$ and, for all $n \in \mathbb{N}$,*

$$(23) \qquad X_n \leq_{\mathcal{D}} \tfrac{1}{\sqrt{2}} \max\{X_{I_n}, X'_{n-1-I_n}\} + A_n,$$

*where $(X_n)_{n \in \mathbb{N}_0}, (X'_n)_{n \in \mathbb{N}_0}, I_n$ are independent, $(X'_n)_{n \in \mathbb{N}_0} =_{\mathcal{D}} (X_n)_{n \in \mathbb{N}_0}$ and $I_n \sim \text{unif}\{0, \ldots, n-1\}$. In this situation, if $(A_n)_{n \in \mathbb{N}_0}$ has uniformly exponentially decreasing tails, then so does $(X_n)_{n \in \mathbb{N}_0}$.*

PROOF.   Suppose that, for some $\kappa > 0$ and $C < \infty$,

$$P(A_n \geq x) \leq C \exp(-\kappa x) \qquad \text{for all } x \geq 0, n \in \mathbb{N}.$$

We have to show that there are finite constants $\tilde{\kappa} > 0$ and $\tilde{C} < \infty$ such that

$$P(X_n \geq x) \leq \tilde{C} \exp(-\tilde{\kappa} x) \qquad \text{for all } x \geq 0, n \in \mathbb{N}_0.$$

If we want to prove this by induction, then the case $n = 0$ is clear as $X_0 \equiv 0$. We may also assume that $x \geq x_0 := \log(\tilde{C})/\tilde{\kappa}$ as otherwise the upper bound is greater than 1. Using (23), we obtain

$$P(X_n \geq x) \leq 2 \max_{k=0,\ldots,n-1} P\left(X_k \geq \frac{3x}{2\sqrt{2}}\right) + P\left(A_n \geq \frac{x}{4}\right),$$

so the induction step will work if $\tilde{\kappa}$ and $\tilde{C}$ can be chosen such that

$$2\tilde{C} \exp\left(-\frac{3\tilde{\kappa}x}{2\sqrt{2}}\right) + C \exp\left(-\frac{\kappa x}{4}\right) \leq \tilde{C} \exp(-\tilde{\kappa} x) \qquad \text{for all } x \geq x_0.$$

This can obviously be done if we first choose $\tilde{\kappa} := \kappa/8$, for example, and then choose $\tilde{C}$ large enough.   □

We next translate the basic recursion (3) for the raw silhouette into a recursion for the sequence $(Z_n)_{n \in \mathbb{N}}$ of processes defined in (15). For this, we require the two linear operators $A, B : C_{00}[0,1] \to C_{00}[0,1]$ given by

$$Af(t) := \tfrac{1}{2}f(2t \wedge 1), \qquad Bf(t) := \tfrac{1}{2}f((2t-1)^+) \qquad \text{for all } t \in [0,1]$$

and the function $\phi : [0,1] \to \mathbb{R}$,

$$\phi(t) := \tfrac{1}{2}(t \wedge (1-t)), \qquad 0 \leq t \leq 1.$$

LEMMA 6. *For all $n \in \mathbb{N}$,*

$$\binom{Y^{\circ}(T_n)}{\eta^{\circ}(T_n)} =_{\mathcal{D}} \begin{pmatrix} AY^{\circ}(T_{I_n}) + BY^{\circ}(T'_{n-1-I_n}) + (\eta^{\circ}(T_{I_n}) - \eta^{\circ}(T'_{n-1-I_n}))\phi \\ + (H(I_n) - H(n-1-I_n))\phi \\ 1 + \frac{1}{2}(\eta^{\circ}(T_{I_n}) + \eta^{\circ}(T_{n-1-I_n})) \\ + \frac{1}{2}(H(I_n) + H(n-1-I_n)) - H_n \end{pmatrix},$$

*where $(T_n)_{n \in \mathbb{N}_0}$, $(T'_n)_{n \in \mathbb{N}_0}$ and $I_n$ are independent, $(T_n)_{n \in \mathbb{N}_0} =_{\mathcal{D}} (T'_n)_{n \in \mathbb{N}_0}$ and $I_n \sim \mathrm{unif}(\{0, \dots, n-1\})$.*

PROOF. For a fixed nonempty tree $T$, the basic recurrence (3) gives

$$Y_t(T) = t + \tfrac{1}{2} Y_{2t \wedge 1}(L(T)) + \tfrac{1}{2} Y_{(2t-1)^+}(R(T))$$

and a straightforward calculation results in

$$Y_t^{\circ}(T) = (AY^{\circ}(L(T)))_t + (BY^{\circ}(R(T)))_t + (\eta(L(T)) - \eta(R(T)))\phi(t).$$

Similarly,

$$\eta^{\circ}(T) = 1 + \tfrac{1}{2}(\eta(L(T)) + \eta(R(T))) - H(\#T).$$

From these, the statement follows on using the distributionally recursive structure of BST sequences explained at the beginning of Section 4. $\square$

We now introduce the space $\mathcal{M}$ of probability measures $\mu$ on the Borel subsets of the space $S$ defined in (16) that satisfy the conditions

$$\int (\|f\|_\infty^2 + x^2) \mu(df, dx) < \infty \quad \text{and} \quad \int x \mu(df, dx) = 0.$$

On $\mathcal{M}$, we define a metric $d$ by

$$d(\mu, \nu)^2 := \inf\{\max\{E\|Y - \bar{Y}\|_\infty^2, 7E(\eta - \bar{\eta})^2\} : (Y, \eta) \sim \mu, (\bar{Y}, \bar{\eta}) \sim \nu\}.$$

The factor 7 will be useful in the proof of Lemma 7 below. As at the end of the proof of Theorem 2, we now construct a (nonlinear) operator $\Psi: \mathcal{M} \to \mathcal{M}$ whose definition is motivated by passing to the limit in the recursion given in Lemma 6: for $\mu \in \mathcal{M}$, let $\Psi(\mu)$ be the joint distribution of the random function

$$AY + BY' + (\eta - \eta')\phi + (\log \xi - \log(1 - \xi))\phi$$

and the real random variable

$$1 + \tfrac{1}{2}(\eta + \eta') + \tfrac{1}{2}(\log \xi + \log(1 - \xi)),$$

where $(Y, \eta)$, $(Y', \eta')$ and $\xi$ are independent, $(Y, \eta) \sim \mu$, $(Y', \eta') \sim \mu$ and $\xi \sim \mathrm{unif}(0, 1)$. It is easy to check that $\Psi$ indeed maps $\mathcal{M}$ into $\mathcal{M}$.

LEMMA 7. *$\Psi$ is a strong contraction on $(\mathcal{M}, d)$.*

Proof.    Let $\mu$ and $\nu$ be elements of $\mathcal{M}$. For any given $\varepsilon > 0$, we can find $(Y, \eta) \sim \mu$ and $(\bar{Y}, \bar{\eta}) \sim \nu$ such that

$$\max\{E\|Y - \bar{Y}\|_\infty^2, 7E(\eta - \bar{\eta})^2\} \leq d(\mu, \nu)^2 + \varepsilon.$$

Now, let $(Y', \eta', \bar{Y}', \bar{\eta}')$ be an independent copy of $(Y, \eta, \bar{Y}, \bar{\eta})$ and let $\xi \sim$ unif$(0, 1)$ be independent of the two random quantities $(Y, \eta, \bar{Y}, \bar{\eta})$ and $(Y', \eta', \bar{Y}', \bar{\eta}')$. By the definition of the operator $\Psi$ and the metric $d$,

$$d(\Psi(\mu), \Psi(\nu))^2$$
$$\leq \max\{E\|A(Y - \bar{Y}) + B(Y' - \bar{Y}') + ((\eta - \bar{\eta}) - (\eta' - \bar{\eta}'))\phi\|_\infty^2,$$
$$7E(\tfrac{1}{2}((\eta - \bar{\eta}) + (\eta' - \bar{\eta}')))^2\}.$$

For the second component, we use independence of $\eta - \bar{\eta}$ and $\eta' - \bar{\eta}'$ and the fact that both have the same distribution to obtain

$$7 \cdot E(\tfrac{1}{2}((\eta - \bar{\eta}) + (\eta' - \bar{\eta}')))^2 = \tfrac{7}{2} \cdot E(\eta - \bar{\eta})^2 \leq \tfrac{1}{2}(d(\mu, \nu)^2 + \varepsilon).$$

The starting point for a similar analysis of the more complicated first part is the observation that $Af$ vanishes on $[\tfrac{1}{2}, 1]$ and that $Bf$ vanishes on $[0, \tfrac{1}{2}]$ for all $f \in C_{00}[0, 1]$ so that, splitting the supremum accordingly,

$$\|A(Y - \bar{Y}) + B(Y' - \bar{Y}') + ((\eta - \bar{\eta}) - (\eta' - \bar{\eta}'))\phi\|_\infty^2$$
$$\leq \max\{\|A(Y - \bar{Y}) + ((\eta - \bar{\eta}) - (\eta' - \bar{\eta}'))\phi\|_\infty^2,$$
$$\|B(Y' - \bar{Y}') + ((\eta - \bar{\eta}) - (\eta' - \bar{\eta}'))\phi\|_\infty^2\}.$$

The sum of the two terms provides an upper bound for the maximum, hence,

$$E\|A(Y - \bar{Y}) + B(Y' - \bar{Y}') + ((\eta - \bar{\eta}) - (\eta' - \bar{\eta}'))\phi\|_\infty^2$$
$$\leq E\|A(Y - \bar{Y}) + ((\eta - \bar{\eta}) - (\eta' - \bar{\eta}'))\phi\|_\infty^2$$
$$+ E\|B(\bar{Y} - \bar{Y}) + ((\eta - \bar{\eta}) - (\eta' - \bar{\eta}'))\phi\|_\infty^2.$$

The two terms on the right-hand side have the same structure. Using the triangle inequality for the supremum norm, $\|\phi\|_\infty \leq \tfrac{1}{4}$, $\|Af\|_\infty \leq \tfrac{1}{2}\|f\|_\infty$ and Minkowski's inequality, we obtain

$$E\|A(Y - \bar{Y}) + ((\eta - \bar{\eta}) - (\eta' - \bar{\eta}'))\phi\|_\infty^2$$
$$\leq E(\tfrac{1}{2}\|Y - Y'\|_\infty + \tfrac{1}{4}|(\eta - \bar{\eta}) - (\eta' - \bar{\eta}')|)^2$$
$$\leq (\tfrac{1}{2}(E\|Y - Y'\|_\infty^2)^{1/2} + \tfrac{1}{4}(E((\eta - \bar{\eta}) - (\eta' - \bar{\eta}'))^2)^{1/2})^2.$$

Inside the large outer brackets, we now use

$$E\|Y - Y'\|_\infty^2 \leq d(\mu, \nu)^2 + \varepsilon$$

and

$$E((\eta - \bar\eta) - (\eta' - \bar\eta'))^2 = 2E(\eta - \bar\eta)^2 \le \tfrac{2}{7}(d(\mu,\nu)^2 + \varepsilon),$$

which, combined, lead to the upper bound

$$E\|A(Y - \bar Y) + ((\eta - \bar\eta) - (\eta' - \bar\eta'))\phi\|_\infty^2$$

$$\le (\tfrac{1}{2}(d(\mu,\nu)^2 + \varepsilon)^{1/2} + \tfrac{1}{\sqrt{56}}(d(\mu,\nu)^2 + \varepsilon)^{1/2})^2$$

$$\le c(d(\mu,\nu)^2 + \varepsilon)$$

with some constant $c < 1/2$. Using the same arguments with the terms involving the operator $B$, we arrive at

$$d(\Psi(\mu), \Psi(\nu))^2 \le 2c(d(\mu,\nu)^2 + \varepsilon).$$

Since $c$ does not depend on $\varepsilon$, we now obtain the strong contraction property on letting $\varepsilon$ tend to 0. $\square$

THEOREM 8. *With $Z_n$, $n \in \mathbb{N}$, and $\Psi$ as above, $Z_n \to_{\mathcal{D}} Z_\infty$ as $n \to \infty$, where the distribution of $Z_\infty$ is the unique fixed point of the operator $\Psi$.*

PROOF. Let

$$W_n(\delta) := \sup_{\substack{0 \le s,t \le 1 \\ |s-t| \le \delta}} |Y_t^\circ(T_n) - Y_s^\circ(T_n)|, \qquad \delta \ge 0,$$

be the modulus of continuity of the process $Y^\circ(T_n)$. Using $|\phi(s) - \phi(t)| \le |s-t|/2$ and Lemma 6, we obtain

$$(24) \quad \begin{aligned} W_n(\delta) \le_{\mathcal{D}} &\tfrac{1}{2}\max\{W_{I_n}(2\delta), W'_{n-1-I_n}(2\delta)\} \\ &+ \tfrac{\delta}{2}|\eta^\circ(T_{I_n}) + \eta^\circ(T_{n-1-I_n})| + \tfrac{\delta}{2}|H(I_n) - H(n-1-I_n)|, \end{aligned}$$

where $(W_n)_{n\in\mathbb{N}}$, $(W'_n)_{n\in\mathbb{N}}$ and $I_n$ are independent, $(W_n)_{n\in\mathbb{N}} =_{\mathcal{D}} (W'_n)_{n\in\mathbb{N}})$ and $I_n \sim \mathrm{unif}(\{0,\dots,n-1\})$. Now, let

$$\tilde W_n := \sup_{0 < \delta \le 1} \delta^{-1/2} W_n(\delta).$$

Clearly, $\tilde W_0 \equiv 0$, and (24) implies that

$$\tilde W_n \le_{\mathcal{D}} \tfrac{1}{\sqrt{2}}\max\{\tilde W_{I_n}, \tilde W'_{n-1-I_n}\} + A_n + B_n \qquad \text{for all } n \in \mathbb{N}$$

with

$$A_n := \tfrac{1}{2}|\eta^\circ(T_{I_n}) + \eta^\circ(T_{n-1-I_n})|, \qquad B_n := |H(I_n) - H(n-1-I_n)|$$

and the usual distributional assumptions. Equation (13) in Theorem 3 implies that $(A_n)_{n\in\mathbb{N}}$ has uniformly exponentially decreasing tails. To obtain the same property for $B_n$, we first observe that it is enough to treat $H(n) - H(I_n)$, which is nonnegative. It is easy to check that $H(n)1_{\{I_n=0\}}$ has uniformly exponentially decreasing tails. Further, we have

$$(H(n) - H(I_n))1_{\{I_n \neq 0\}} \leq_{\mathcal{D}} -\log(\xi) + 1$$

with $\xi \sim \text{unif}(0, 1)$, so this term also has the required tail property. Therefore, Lemma 5 can be applied, leading to a uniform upper bound for the tails of $\tilde{W}_n$. In particular,

$$\limsup_{\delta \downarrow 0} \; W_n(\delta) = 0 \qquad \text{in probability,}$$

which shows that the sequence $(Y^\circ(T_n))_{n \in \mathbb{N}}$ is tight in $C_{00}[0, 1]$; see Section 8 in Billingsley (1968). In view of Theorem 2, $(\eta^\circ(T_n))_{n \in \mathbb{N}}$ is also tight in $C_{00}[0, 1]$, which, by a standard argument using Prohorov's theorem, implies tightness of the sequence $(Z_n)_{n \in \mathbb{N}}$. Convergence of the finite-dimensional distribution was obtained in Theorem 4. Combining elements, we see that $Z_n \to_{\mathcal{D}} Z_\infty$ for some $S$-valued random quantity $Z_\infty$. Finally, we can pass to the limit in the distributional equation given in Lemma 6 and then use Lemma 7, together with Banach's fixed point theorem to complete the proof. $\square$

What can be said about the paths of the process part $Y_\infty$ of the limit $Z_\infty$? The functional limit theorem implies that the paths are continuous as everything happens in $C([0, 1])$. The maximum of the raw silhouette is the height $H_n$ of the tree. It has been shown by Devroye (1986) that $H_n / \log n$ converges to $c_+ = 4.311\ldots$ in probability as $n \to \infty$. Similarly, the minimum is the tree's fill (or saturation) level $L_n$ and $L_n / \log_n$ converges to $c_- = 0.373\ldots$ almost surely as $n \to \infty$, a result obtained by Biggins (1997) in the context of branching random walks. Heuristically, as $c_- < c_+$, we would therefore expect that the paths of the limit process are not differentiable. As the next result shows, almost all paths of $Y_\infty$ are indeed not even Lipschitz, but they are Hölder continuous of order $\alpha$ for all $\alpha < 1$.

THEOREM 9. *With* $Z_\infty = (Y_\infty, \zeta)$ *as in Theorem* 8,

$$P\left(\sup_{0 \leq s < t \leq 1} \frac{|Y_\infty(t) - Y_\infty(s)|}{|t - s|} = \infty\right) = 1$$

*and, for all* $\alpha < 1$,

$$P\left(\sup_{0 \leq s < t \leq 1} \frac{|Y_\infty(t) - Y_\infty(s)|}{|t - s|^\alpha} < \infty\right) = 1.$$

PROOF. Taking $s = 0$ and $t = 2^{-k}$ in (17), we obtain

$$(25) \qquad \frac{|Y_\infty(t) - Y_\infty(s)|}{|t - s|} =_{\mathcal{D}} |k + \rho_k + \eta|,$$

where $-\rho_k$ has a gamma distribution with shape parameter $k$ and scale parameter 1, and $\eta$ is random variable whose distribution does not depend on $k$. For $k$ large, $(k + \rho_k)/\sqrt{k}$ is close in distribution to a standard normal, which shows that the right-hand side of (25) tends to $\infty$ in probability as $k \to \infty$. This proves the first statement.

For proof of the second part, we use the Kolmogorov–Chentsov theorem; see, for example, Kallenberg (1997), Theorem 2.23. Let $\rho_k$ and $\eta$ be as above; note that $-\rho_k$ can be written as the sum of $k$ independent random variables that are exponentially distributed with parameter 1. Moments of all orders exist for these and for $\zeta$. Hence, by Minkowski's inequality, for each $l \in \mathbb{N}$, there is a constant $C_l$ such that

$$E|k + \rho_k + \eta|^l \leq C_l k^l \qquad \text{for all } k \in \mathbb{N}.$$

If $s$ and $t$ are such that $s = j2^{-k}$ and $t = (j+1)2^{-k}$ for some $j \in \{0, \ldots, 2^k - 1\}$, then we obtain, again using (17),

$$(26) \qquad E|Y_\infty(t) - Y_\infty(s)|^l \leq |t - s|^l |\log_2 |t - s||^l C_l.$$

The desired Hölder continuity now follows on using two facts: first, the right-hand side of (26) can be bounded by $|t - s|^{l+\delta}$ for all $\delta > 0$ and we may choose an arbitrarily large value for $l$; second, the chaining proof of the Hölder part of the Kolmogorov–Chentsov theorem makes use of the moment bounds only for values of $s$ and $t$ that are of the above form, that is, binary rational neighbors. $\square$

**6. Remarks.** In this final section, we briefly discuss another family of search trees, comment on the methodology and close with a final remark on the "big picture."

(i) As in the previous sections, let $(\xi_n)_{n \in \mathbb{N}}$ be a sequence of independent random variables, all uniformly distributed on the unit interval. The DST (digital search tree) algorithm uses the binary expansion of the values as a directive of how to travel through the binary tree, storing each value in the first free (i.e., the unique external) node; again, we refer to Knuth (1973), Mahmoud (1992) and Sedgewick and Flajolet (1996) for more information. As in the BST case, the algorithm produces a sequence $(T_n)_{n \in \mathbb{N}}$ of random trees, where $T_n$ is the DST output for the first $n$ variables $\xi_1, \ldots, \xi_n$ of the sequence. In contrast to the BST situation, we no longer have invariance of the resulting random structures under strictly monotone transformations of the input data. However, we still have a simple "stochastic dynamics": in both cases, $T_{n+1}$ is obtained by adding a randomly selected element of $\partial T_n$ to $T_n$, but, whereas in the BST case, one of the $n+1$ external nodes of $T_n$ is chosen uniformly, in the DST case, it is chosen with probability $2^{-k}$, where $k$ is the height of the external node (the fact that these values sum to 1 has already been mentioned in the proof of Theorem 2). The silhouette of such a tree, in raw and in normalized integrated form, is shown in Figure 2. It is "visually obvious" when comparing this with Figure 1 that these functions are quite different in the BST and DST cases (Figures 1 and 2 are based on the same input and use the same scale).

A first point of interest is the fact that the associated discounted external path length $\eta_n$ now has a direct algorithmic interpretation as the conditional expectation of the number of bit checks necessary to insert the next data value $\xi_{n+1}$ into the tree $T_n$, given $T_n$. Generally, the DST output is closer to the "ideal" search tree of minimal height, but makes stronger assumptions on the nature of the input. One indication of this is the fact that distributional fluctuations appear in the DST situation; indeed, if we always stored the next item in an external node of minimal height, we would have

$$\eta_{n,\text{opt}} - \log_2 n = \phi(\{\log_2 n\}),$$

where $\{x\}$ denotes the fractional part of $x$ and $\phi(x) := 2^x - 1 - x$, $0 \le x \le 1$. Some heuristic arguments support the conjecture that the expectation of $\eta_n$ differs from $\eta_{n,\text{opt}}$ by an asymptotically negligible amount and that the variance of $\eta_n$ tends to 0 as $n \to \infty$.

As a second point of interest, we note that the analysis of the associated silhouette processes begins to bifurcate at the earliest possible point, that is, in the situation considered in Theorem 1. As a result of the stochastic dynamics explained above, the movement along a particular path $s \in [0,1)$, that is, the behaviour of $X_s(T_n)$ for $n = 1, 2, \ldots$, is a Markov chain of pure-birth type with state space $\mathbb{N}$ and birth rates $p_{k,k+1} = 2^{-k}$. The associated distributions converge along suitably chosen subsequences if we simply subtract $\log n$; see Dennert and Grübel (2007) for a recent probabilistic approach. Recall that in the BST case, the variance of $X_s(T_n)$ grows at a logarithmic rate and that a suitably rescaled version of $X_s(T_n)$ is asymptotically normal. Moreover, the random variables $X_s(T_n)$ and $X_t(T_n)$ are asymptotically independent by Theorem 1 if $s \neq t$; in particular, the absolute difference between the two converges to $\infty$ in probability as $n \to \infty$. In the DST case, however, it follows easily from the result mentioned above, that the family of distributions of the differences $X_t(T_n) - X_s(T_n)$, $n \in \mathbb{N}$, is tight. Again on the basis of heuristic arguments, I conjecture that the distributional periodicities disappear in an appropriately standardized version of the silhouette, such as $(X_t(T_n) - X_0(T_n))_{0 \leq t \leq 1}$.

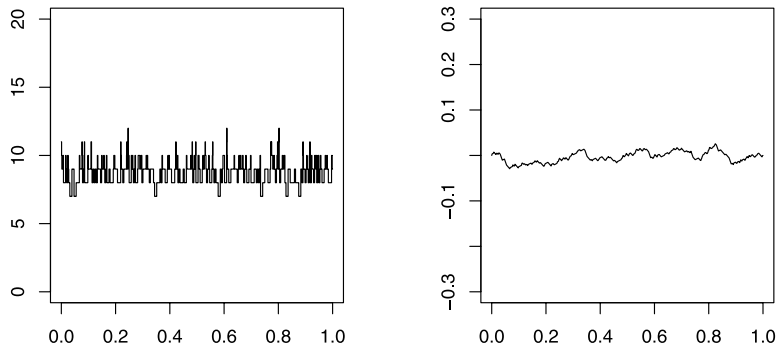

FIG. 2.   *An example of a silhouette (left) and the corresponding normalized integrated silhouette (right) in the DST case.*

(ii) In our proofs, we have used martingale results and contraction arguments. A survey of the contraction method in the context of the analysis of algorithms is given in Rösler and Rüschendorf (2001) and, with emphasis on the multivariate case, in Neininger and Rüschendorf (2006). A first use of the contraction method in connection with the analysis of algorithms on the level of stochastic processes, as in the present paper, can be found in Grübel and Rösler (1996). Roughly, martingale arguments often provide almost sure convergence in cases where the contraction method only yields convergence in distribution, but the latter seems to have advantages if our interest is in the properties of the limit distribution. The two methods are closely related to the complementary aspects of BST sequences, the dynamic structure and the distributionally recursive structure, that we have used repeatedly in the previous sections.

(iii) According to Knuth (1997), page 308, trees are "the most important nonlinear structures that arise in computer algorithms." Given a sequence of input data, both the BST and the DST algorithms generate a sequence $(T_n)_{n \in \mathbb{N}}$ of binary trees that grow by one node at a time. As has been shown in Luczak and Winkler (2004), even in the case of uniformly distributed trees, there is a dynamical procedure that builds these structures in this sequential manner. From that point of view, in all three cases, the stochastic process $(T_n)_{n \in \mathbb{N}}$ of trees is a transient Markov chain with a denumerable state space $E$, with $E$ the set of all finite and prefix-stable subsets of the denumerable set $\mathcal{N}$ of nodes. One would expect that, in a rough sense, the limit is always the complete binary tree $T_\infty$. However, this is not true in the uniform case; see Luczak and Winkler (2004) and the references given therein. However, for the search trees that we have considered in the present paper, the fill level converges to $\infty$ with probability 1, so a simple compactification $E_\infty := E \cup \{T_\infty\}$ makes $(T_n)_{n \in \mathbb{N}}$ a sequence that converges with probability 1 – if convergence means that every $u \in \mathcal{N}$ will eventually be an element of $T_n$. From a general theoretical point of view, the results of the present paper can be regarded as a first step toward a more detailed asymptotic analysis, going beyond the one-point compactification toward classical notions such as the Martin boundary; see Sawyer (1997). In this connection, it is interesting to note that Régnier's (1989) analysis of Quicksort is based on one specific harmonic function associated with the Markov chain that arises in the BST case.

## REFERENCES

ALDOUS, D. (1991). The continuum random tree. II. An overview. In *Stochastic Analysis (Durham, 1990). London Mathematical Society Lecture Note Series* **167** 23–70. Cambridge Univ. Press, Cambridge. MR1166406

BIGGINS, J. D. (1997). How fast does a general branching random walk spread? In *Classical and Modern Branching Processes (Minneapolis, MN, 1994). IMA Vol. Math. Appl.* **84** 19–39. Springer, New York. MR1601689

BILLINGSLEY, P. (1968). *Convergence of Probability Measures.* Wiley, New York. MR0233396

CHAUVIN, B., DRMOTA, M. and JABBOUR-HATTAB, J. (2001). The profile of binary search trees. *Ann. Appl. Probab.* **11** 1042–1062. MR1878289

CRAMER, M. (1996). A note concerning the limit distribution of the quicksort algorithm. *RAIRO Inform. Théor. Appl.* **30** 195–207. MR1415828

DENNERT, F. and GRÜBEL, R. (2007). Renewals for exponentially increasing lifetimes, with an application to digital search trees. *Ann. Appl. Probab.* **17** 676–687. MR2308339

DEVROYE, L. (1986). A note on the height of binary search trees. *J. Assoc. Comput. Mach.* **33** 489–498. MR849025

DEVROYE, L., FILL, J. A. and NEININGER, R. (2000). Perfect simulation from the Quicksort limit distribution. *Electron. Comm. Probab.* **5** 95–99 (electronic). MR1781844

FILL, J. A. and JANSON, S. (2000). Smoothness and decay properties of the limiting Quicksort density function. In *Mathematics and Computer Science (Versailles, 2000)* 53–64. Birkhäuser, Basel. MR1798287

GRÜBEL, R. (2005). A hooray for Poisson approximation. In *2005 International Conference on Analysis of Algorithms. Discrete Math. Theor. Comput. Sci. Proc.* **AD** 181–191 (electronic). Assoc. Discrete Math. Theor. Comput. Sci., Nancy. MR2193118

GRÜBEL, R. and RÖSLER, U. (1996). Asymptotic distribution theory for Hoare's selection algorithm. *Adv. in Applied Prob.* **28** 252–269.

KALLENBERG, O. (1997). *Foundations of Modern Probability.* Springer, New York. MR1464694

KNUTH, D. E. (1973). *The Art of Computer Programming. Volume 3.* Addison-Wesley, Reading, MA. MR0445948

KNUTH, D. E. (1975). *The Art of Computer Programming, Volume 1*, 3nd ed. Addison-Wesley, Reading, MA.

LE GALL, J.-F. (1999). *Spatial Branching Processes, Random Snakes and Partial Differential Equations.* Birkhäuser, Basel. MR1714707

LUCZAK, M. and WINKLER, P. (2004). Building uniformly random subtrees. *Random Structures Algorithms* **24** 420–443. MR2060629

MAHMOUD, H. M. (1992). *Evolution of Random Search Trees.* Wiley, New York. MR1140708

NEININGER, R. and RÜSCHENDORF, L. (2006). A survey of multivariate aspects of the contraction method. *Discrete Math. Theor. Comput. Sci.* **8** 31–56 (electronic). MR2247515

PITMAN, J. (2006). *Combinatorial Stochastic Processes. Lecture Notes in Math.* **1875**. Springer, Berlin. MR2245368

PITTEL, B. (1985). Asymptotical growth of a class of random trees. *Ann. Probab.* **13** 414–427. MR781414

PITTEL, B. (1986). Paths in a random digital tree: Limiting distributions. *Adv. in Appl. Probab.* **18** 139–155. MR827333

RÉGNIER, M. (1989). A limiting distribution for quicksort. *RAIRO Inform. Théor. Appl.* **23** 335–343. MR1020478

RÖSLER, U. (1991). A limit theorem for "Quicksort". *RAIRO Inform. Théor. Appl.* **25** 85–100. MR1104413

RÖSLER, U. and RÜSCHENDORF, L. (2001). The contraction method for recursive algorithms. *Algorithmica* **29** 3–33. Average-case analysis of algorithms (Princeton, NJ, 1998). MR1887296

SAWYER, S. A. (1997). Martin boundaries and random walks. *Contemp. Math.* **206** 17–44. MR1463727

SEDGEWICK, R. and FLAJOLET, P. (1996). *An Introduction to the Analysis of Algorithms.* Addison-Wesley, Reading, MA.

INSTITUT FÜR MATHEMATISCHE STOCHASTIK
LEIBNIZ UNIVERSITÄT HANNOVER
POSTFACH 60 09
D-30060 HANNOVER
GERMANY
E-MAIL: rgrubel@stochastik.uni-hannover.de